\documentclass[10pt,psamsfonts]{amsart}

\newtheorem{theorem}{Theorem}[section]
\newtheorem{proposition}[theorem]{Proposition}
\newtheorem{lemma}[theorem]{Lemma}

\theoremstyle{definition}
\newtheorem{definition}[theorem]{Definition}
\newtheorem{remark}[theorem]{Remark}

\newcommand{\R}{\mathbb{R}}
\newcommand{\N}{\mathbb{N}}
\newcommand{\ep}{\varepsilon}

\numberwithin{equation}{section}

\def\porz(#1){\gamma\langle#1\rangle}
\def\uporz(#1){\overline{p}\langle#1\rangle}
\def\lporz(#1){\underline{p}\langle#1\rangle}

\def\por(#1){\gamma(#1)}
\def\upor(#1){\overline{p}(#1)}
\def\lpor(#1){\underline{p}(#1)}

\def\porm(#1){\gamma(#1)}
\def\uporm(#1){\overline{\mathrm{por}}(#1)}
\def\lporm(#1){\underline{\mathrm{por}}(#1)}
\def \cal{\mathcal}
\def \dist{\mathrm{dist}}

\def \span{\mathrm{sp}}
\def \vn{\mathrm{int}}
\def \wt{\widetilde}
\def \N{{\bf N}}
\def \R{{\bf R}}
\def \ep{\varepsilon}
\def \vf{\varphi}
\def \DC{\mathcal{DC}}
\def \L{\mathcal{L}}
\def \A{\mathcal{A}}

\begin{document}

\renewcommand{\labelenumi}{{\rm(\roman{enumi})}}


\def\MR{{\bf MR\/}~}

\centerline{ON LIPSCHITZ AND D.C. SURFACES OF FINITE CODIMENSION IN A BANACH SPACE}

\bigskip
\bigskip

\centerline{LUD\v EK ZAJ\' I\v CEK, Praha}



\bigskip
\bigskip

{\it Abstract.}\  Properties of Lipschitz and d.c. surfaces of finite codimension in a 
Banach space, and properties of generated $\sigma$-ideals
 are studied. These $\sigma$-ideals naturally appear in the differentiation theory and in the abstract approximation theory. Using these properties, we improve an unpublished 
  result of M. Heisler which gives an alternative proof of a result 
  of D. Preiss on singular points of convex functions.

\bigskip
\bigskip
 
 {\it Keywords}:\  Banach space, Lipschitz surface, d.c. 
surface, multiplicity points  of monotone operators, singular points of convex functions, Aronszajn null sets.

\bigskip

{\it  MSC 2000}:\ 46T05, 58C20, 47H05
\bigskip
\bigskip

\section{Introduction}

Let $X$ be a real separable Banach space. A number of 
$\sigma$-ideals of subsets of $X$ are considered in the literature. Besides the most classical
  system of first category sets mention the $\sigma$-ideals of
   Haar null sets, Aronszajn (equivalently Gaussian) null sets (see \cite{BL}),
$\Gamma$-null sets (see \cite{LP}, \cite{LP01}) and 
 $\sigma$-(lower or upper) porous sets (see e.g., \cite{Zasu}).
In some questions   of the differentiability theory and of the abstract 
approximation theory,  the $\sigma$-ideals  $\L^1(X)$ and $\DC^1(X)$
 generated by Lipschitz and d.c. Lipschitz hypersurfaces (i.e., ``graphs''
 of Lipschitz and of d.c. Lipschitz functions), respectively, naturally
  appear.  These $\sigma$-ideals are proper subsystems of all 
  $\sigma$-ideals mentioned above.
  The sets from $\L^1(X)$ 
 were used in $\R^2$ (under a different but equivalent definition) by
 W.H. Young (under the name ``ensemble rid\' ee'') and by H. Blumberg
(under the name ``sparse set''); cf. \cite[p. 294]{Za83}.
 These sets were used in $\R^n$ e.g., (implicitly) by P. Erd\" os \cite{Er},
 and in infinite-dimensional spaces (possibly for the first time)
  in \cite{Za78b} and \cite{Za78a}.
The sets from $\DC^1(X)$ were probably first applied in  \cite{Za79} (cf. \cite[p. 93]{BL}).
In some articles (e.g., \cite{Za78b}, \cite{Za79}, \cite{Za83}, \cite{Ve86})
 also sets from smaller $\sigma$-ideals $\L^n(X)$ and $\DC^n(X)$
 generated by Lipschitz and d.c. Lipschitz surfaces of codimension
  $n>1$ were used.

In the present article we prove some properties of Lipschitz and Lipschitz locally d.c. surfaces
 of finite codimension  (Section 3; Proposition \ref{brower} and Proposition \ref{per}).
 
  Using these properties, we study in Section 4
  sets which are projections of sets from $\L^n(X)$ on a closed space $Y \subset X$ of codimension
   $d<n$. The study of such projections was suggested by D. Preiss in connection with a result of \cite{PreH} (see Remark \ref{DH}(i)). M. Heisler \cite{HeD} proved that any such projection
     is a first category set in $Y$, which provides (together with a result of \cite{Za79}) an alternative proof
      of a result of \cite{PreH}. We prove that each such projection is also a subset
       of an Aronszajn null set in $Y$ (and even a subset of a set from a smaller class 
        $\cal C_{n}^*$). As a consequence, we obtain  a result on projections of sets of multiplicity
          of monotone operators (Theorem \ref{mon}) which  improves both \cite[Theorem 1.3.]{PreH} and the corresponding result of \cite{HeD}.
        
        Our proof is more transparent than that of \cite{HeD} and gives stronger results, since it
         uses ``perturbation'' Proposition  \ref{per}. To prove (and apply) it, we need some results on perturbations
          of finite-dimensional subspaces. These results are collected in Preliminaries, where
           also needful results on d.c. mappings are  recalled.

\section{Preliminaries}

We consider only real Banach spaces. By $\span\{M\}$ we denote the linear span of the set $M$.
A mapping is called $K$-Lipschitz if it is Lipschitz with a (not necessarily minimal) constant $K$. 
 A bijection $f$ is called bilipschitz ($K$-bilipschitz) if both $f$ and $f^{-1}$ are Lipschitz ($K$-Lipschitz).

A real function on an open convex subset of a Banach space  is called d.c. (delta-convex) if it is a difference of two continuous convex functions. Hartman's notion of d.c. mappings between Euclidean spaces \cite{Ha} was generalized
 and studied in \cite{VeZa}. 
 \begin{definition}\label{dc}
 Let $X, Y$ be Banach spaces, $C\subset X$  an open convex set, and let $F: C \to Y$ be a continuous mapping. We say that
  $F$ is d.c. if there exists a continuous convex function $f: C \to \R$ such that $y^* \circ F + f$ is convex whenever  $y^* \in Y^*$, $\|y^*\| \leq 1$. 
 \end{definition}
It is easy to see (cf. \cite[Corollary 1.8.]{VeZa}) that, if $Y$ is finite dimensional, then $F$ is d.c. if and only if $y^* \circ F$ is d.c. for each
 $y^* \in Y^*$ (or for each $y^*$ from a fixed basis of $Y^*$). Note also that each d.c. mapping is locally
  Lipschitz (\cite[p. 10]{VeZa}). If $X$ is finite-dimensional, then each locally d.c. mapping is d.c. (see. \cite[p. 14]{VeZa}) but it is not true
   (see \cite{KM}) if $X$ is infinite dimensional. We will need also the following well-known facts on d.c. mappings.
   \begin{lemma}\label{fakta}
   Let $X$, $X_1$, $Y$, $Y_1$, $Y_2$, $Z$ be Banach spaces. 
   \begin{enumerate}
   \item
   Let $f: X \to Y$ be d.c. and  let $g: X_1 \to X$, $h: Y \to Y_1$ be linear and continuous. Then
    both $f\circ g$ and $h \circ f$ are d.c.
   \item
   A mapping $f= (f_1,f_2): X \to Y_1 \times Y_2$ is d.c. if and only if both $f_1$ and $f_2$ are d.c.
   \item
   If $g: X \to Y$, $h: X \to Y$ are d.c. and $a,b \in \R$, then $ag + bh$ is d.c.
   \item
   If $f: X \to Y$ is locally d.c. and $g:Y \to Z$ is locally d.c., then $g \circ f$ is locally d.c.
   \item
   Suppose that $G: X \to Y$ is a linear isomorphism, $g: X \to Y$ is a locally d.c.  bilipschitz bijection,
    and the range of $g-G$ is contained in a finite dimensional space. Then $g^{-1}$ is locally d.c.
    \end{enumerate}
   \end{lemma}
\begin{proof}
The statements (i) and (ii) are very easy (cf. \cite[Lemma 1.5. and Lemma 1.7.]{VeZa}) and (iii) follows from
 (i) and (ii). The statement (iv) is a special case of
 \cite[Theorem 4.2.]{VeZa} and (v) is a special case of  \cite[Theorem 2.1.]{D}.
\end{proof}

We will need some notions and results concerning distances of two subspaces of a Banach space, which
 are well-known from the perturbation theory of linear operators (\cite{GK}, \cite{K}, \cite{B}).
 Let $X$ be a Banach space and  $S(X)$ be the unit sphere of $X$. Let $Y$ and $Z$ be closed non-trivial subspaces of $X$. Then the {\it gap
  between $Y$ and $Z$} (called also the opening or the deviation of $Y$ and $Z$) is defined by
  $$    \gamma(Y,Z) = \max \{ \sup_{y \in Y \cap S(X)} \dist (y,Z), \sup_{z \in Z \cap S(X)} \dist (z,Y) \}.$$ 
  We set $\gamma(\{0\},\{0\}) := 0$ and $\gamma(Y,Z)=1$ if one and only one of $Y$, $Z$ is $\{0\}$.
The gap need not be a metric on the set of all non-trivial subspaces of $X$; this property
 has the distance $\rho(Y,Z)$ between $Y$ and $Z$ defined  as the Hausdorff distance between
  $Y \cap S(X)$ and $Z \cap S(X)$.
  
  We will work with the gap $\gamma (Y,Z)$. However, since it is easy to prove (see e.g., \cite{K}) that 
   (for nontrivial $Y$, $Z$) always
  \begin{equation}\label{vztah}
    \rho(Y,Z)/2 \leq   \gamma(Y,Z) \leq \rho(Y,Z),
  \end{equation}
   we could work also with $\rho(Y,Z)$. We will need the following well-known facts.
   
   \begin{lemma}\label{zname}
   Let $X$ be a Banach space and $F$, $\wt F$, $K$ be   finite dimensional subspaces of $X$. Then:
   \begin{enumerate}
   \item[(i)]
   If  $\gamma(F, \wt F) < 1$, then $\dim F = \dim \wt F$.
   \item[(ii)]
   If  $F \cap K = \{0\}$, then there exists $\omega>0$ such that 
    $\gamma(F, \wt F) < \omega$ implies $\wt F \cap K = \{0\}$.
    \item[(iii)]
    If $E \oplus F = X$, then there exists $\omega>0$ such that 
    $\gamma(F, \wt F) < \omega$ implies  $E \oplus \wt F = X$.
  \end{enumerate}
  \end{lemma}
  \begin{proof}
  Statement (i) is proved in \cite{GK} (see \cite[Theorem 2.1]{B}) and (ii) is an easy consequence of
   \eqref{vztah}. (We can also apply \cite[Theorem 5.2]{B} with $Y:= F$, $Z:= K$ and $X:= F \oplus K$.)
    The statement (iii) immediately follows from \cite[Theorem 5.2]{B}.
  \end{proof}

   The following simple lemma is also essentially  well-known. Although it is not stated explicitly in \cite{L}, it follows from \cite[Theorem 2.2.]{L} which works with complex
    Banach spaces. Since the formulation of   \cite[Theorem 2.2.]{L} is rather complicated and we work with real spaces, for the  sake of completness we give a proof.       
   
  \begin{lemma}\label{baze}
  Let $X$ be a Banach space and $(v_1,\dots,v_n)$ be a basis of a space $V \subset X$ and $\ep>0$.  Then there exists $\delta >0$ such that
   the inequalities  $\|w_1-v_1\| < \delta,\dots,  \|w_n-v_n\| < \delta$  imply that $W:= \span\{w_1,\dots,w_n\}$
    is $n$-dimensional  and $\gamma(V,W) < \ep$.          
  \end{lemma} 
  
\begin{proof}
First we will show that there exists $\eta>0$ and $\delta^*>0$ such that the inequality  
\begin{equation}\label{ci}
\left \| \sum_{i=1}^n c_i w_i \right\| \geq \eta  \|c\|_{\infty}
\end{equation}
holds  whenever $\|w_1-v_1\| < \delta^*,\dots,  \|w_n-v_n\| < \delta^*$ and $c=(c_1,\dots,c_n) \in \R^n$ is arbitrary.
To this end observe that there exists $\eta^* >0$ such that \eqref{ci} holds for $\eta=\eta^*$, $w_i=v_i$, and arbitrary $c$.
Put  $\eta := \eta^*/2$ and $\delta^* : = \eta^*/2n$. Then the inequalities  
 $\|w_1-v_1\| < \delta^*,\dots,  \|w_n-v_n\| < \delta^*$  imply that, for each $0 \neq c \in \R^n$, 
 $$ \left \| \sum_{i=1}^n  \frac{c_i}{\|c\|_{\infty}}\, v_i \right\| -
 \left \| \sum_{i=1}^n  \frac{c_i}{\|c\|_{\infty}}\, w_i \right\| \leq    
 \left \| \sum_{i=1}^n  \frac{c_i}{\|c\|_{\infty}}\, (v_i- w_i) \right\| < n \delta^*  = \eta^*/2.        
  $$
Consequently, using the definition of $\eta^*$, we obtain
$$ \left \| \sum_{i=1}^n  \frac{c_i}{\|c\|_{\infty}}\, w_i \right\| \geq
  \left \| \sum_{i=1}^n  \frac{c_i}{\|c\|_{\infty}}\, v_i \right\| - \eta^*/2 \geq \eta^*- \eta^*/2 = \eta,$$
  which implies \eqref{ci}.
  
Now set $\delta:= \min\{\delta^*, \ep \eta/2n\}$ and suppose that the inequalities  
$\|w_1-v_1\| < \delta,\dots,  \|w_n-v_n\| < \delta$ hold.                         
 Let  $w = \sum_{i=1}^n c_i w_i$ with $\|w\|=1$ be given. Set $v = \sum_{i=1}^n c_i v_i$. 
 Since  $\|c\|_{\infty} \leq 1/\eta$ by \eqref{ci}, we obtain
 $$ \|v-w\| \leq \sum_{i=1}^n |c_i| \delta \leq  n (1/\eta) \delta \leq  \ep/2.$$
 Consequently, $\sup_{w \in W \cap S(X)} \dist (w,V) < \ep$. By a quite symmetrical way we obtain        
 $\sup_{v \in V \cap S(X)} \dist (v,W) < \ep$, so
  $\gamma(V,W) < \ep$. Since we can suppose $\ep<1$, we know that $W$ is $n$-dimensional by Lemma \ref{zname}(i).
\end{proof}

\begin{lemma}\label{gaptop}
Let $X$, $Y$ be Banach spaces and $F: X \to Y$ be a linear isomorphism. Then there exists $C>0$ such that
$$ C^{-1} \gamma(F(V),F(W)) \leq \gamma(V,W) \leq C \gamma(F(V), F(W)),$$
whenever $V$ and $W$ are subspaces of $X$.
\end{lemma}
\begin{proof}
We can clearly suppose that $V$ and $W$ are non-trivial. 
Since $F^{-1}$ is also a linear isomorphism, it is clearly sufficient to find $D>0$ such that
$\gamma(V,W) \leq D \gamma(F(V), F(W))$ always holds. Choose $K>0$ such that $F$ is $K$-bilipschitz
 and consider $v \in V$ with  $\|v\|=1$. We can clearly find $\wt  w \in F(W)$ for which
  $\|\wt  w - \|F(v)\|^{-1}\cdot F(v) \| \leq 2 \gamma(F(V), F(W))$. Since
   $\|F(v)\| \leq K$, we have $\| F(v) - \|F(v)\|\cdot \wt  w\| \leq 2K \gamma(F(V), F(W))$, and therefore
    $\|v - \|F(v)\|\cdot F^{-1}(\wt  w)\| \leq 2 K^2 \gamma(F(V), F(W))$. Since the roles of $V$ and $W$ are symmetric,  we can clearly set $D:= 2K^2$.
\end{proof}

\begin{lemma}\label{prpr}
Let $X$ be an infinite dimensional Banach space, $V, W \subset X$  non-trivial finite dimensional spaces, and
  $\delta > 0$. Then there exists a space $\wt V \subset X$ with $\gamma(V, \wt V) < \delta$
  and $\wt V \cap W = \{0\}$.
\end{lemma}
\begin{proof}
Denote $n:= \dim V$, choose an $n$-dimensional space $Y \subset X$ with $Y \cap (V+W) = \{0\}$ and a linear bijection
 $L: V \to Y$. For $t>0$, set $\wt V_t := \{ v + t L(v):\ v \in V\}$. It is easy to check that
  each $\wt V_t$ is an $n$-dimensional space with $\wt V_t \cap W = \{0\}$. 
 Applying Lemma \ref{baze} to a basis $v_1,\dots, v_k$ of $V$ and $w_i := v_i + tL(v_i)$, it is easy
  to see that $\gamma(V, \wt V_t) \to 0$ ($t \to 0+$), which implies our assertion.
\end{proof}

\begin{lemma}\label{kdk}
Let $X$ be a  Banach space, $1 \leq n < \dim X$, and $K \geq 1$.
 Let $X = E \oplus F$, where $F$ is an $n$-dimensional space. Suppose
  that the canonical mapping $\mu: E \oplus F \to E \times F$
   (where $E \times F$ is equipped with the maximum norm) is 
   $K$-bilipschitz. Then there exists $\omega >0$ such that if
 $\wt{F} \subset X$ is a closed space with $\gamma(F,\wt F)< \omega$,
  then $X = E \oplus \wt F$ and  
the canonical mapping $\wt{\mu}: E \oplus \wt 
F \to E \times \wt F$ is $2K$-bilipschitz.
\end{lemma}
\begin{proof}
Distinguishing the cases $\lambda <1$, $\lambda=1$ and
 $\lambda > 1$, it is easy to check that
there exists $1>\omega_0 > 0$  such
  the inequalities
\begin{multline}\label{dve}
K \max(1+\omega, \lambda) + \omega \leq 2K \max(1,\lambda), \\ 
K^{-1} \max(1-\omega, \lambda) - \omega \geq (2K)^{-1} \max(1,\lambda)
\end{multline}
hold for each $\lambda \geq 0$ and $0 < \omega < \omega_0$. 
 By Lemma \ref{zname}(iii), we can choose  $0< \omega < \omega_0$ such that     
   $X = E \oplus \wt F$ whenever $\gamma(F,\wt F)< \omega$.   
Let $\wt F$ with $\gamma(F,\wt F)< \omega$ be given,
and consider arbitrary $\wt f \in \wt F$ and $e \in E$.
 We will prove 
 \begin{equation}\label{bil}
 (2K)^{-1} \max(\|\wt f\|, \|e\|)
  \leq \|\wt f + e\| \leq 2K \max(\|\wt f\|, \|e\|).
\end{equation}
Since the case $\wt f = 0$ is trivial, by homogenuity of
 the norm we can suppose $\|\wt f\| = 1$ and find $f \in F$
  with $\|f - \wt f\| < \omega$. 
   Applying (\ref{dve}) to $\lambda := \|e\|$, we obtain
 \begin{multline*}
 \|\wt f + e\| \leq \|f + e\| + \omega \leq K \max(\|f\|,\|e\|) + 
\omega  \\
\leq  K \max(1+\omega,  \|e\|) + \omega  \leq
 2K \max(1,\|e\|)
 \end{multline*}
and
 $$
 \|\wt f + e\| \geq \|f + e\| - \omega \\ \geq
  K^{-1} \max(1-\omega,\|e\|) -
\omega 
 \geq
 (2K)^{-1} \max(1,\|e\|).
 $$

 Thus,  (\ref{bil}) holds, and $\wt \mu$ is $(2K)$-bilipschitz. 
\end{proof}

\section{Properties of Lipschitz surfaces  of finite codimension}

If $X$ is a Banach space and $X = E \oplus F$, then we denote by
 $\pi_{E,F}$ the projection of $X$ on $E$ along the space $F$.

\begin{definition}
Let $X$ be a Banach space and $A \subset X$.
\begin{enumerate}
\item[(i)]
Let $F$ be a closed subspace of $X$. We say that $A$ is an {\it $F$-Lipschitz
 surface} if there exists a topological complement $E$ of $F$ 
 and a Lipschitz mapping $\vf: E \to F$ such that  $A = \{x+\vf(x): \ x \in E\}$.
 \item[(ii)]
 Let  $1 \leq n < \dim X$ be a natural number. We say that
  $A$ is a {\it Lipschitz surface 
 of codimension $n$} if $A$ is an $F$-Lipschitz surface for some $n$-dimensional  space
 $F \subset X$.
 \item[(iii)]
 If we consider in (i) mappings $\vf: E \to F$ which are d.c. (resp. Lipschitz d.c., locally d.c., Lipschitz locally d.c.), we
  obtain the notions of an $F$-d.c. surface, d.c.  surface 
 of codimension $n$ (resp. $F$-Lipschitz d.c. surface, etc.).
 A Lipschitz surface (resp. d.c. surface, etc.) of codimension 1 is said to be a
{\it Lipschitz hypersurface (resp. d.c. hypersurface, etc.)}.
 \item[(iv)]
The
 $\sigma$-ideals of sets which can be covered by countably many
 Lipschitz surfaces (d.c. surfaces) of codimension $n$ will be denoted
 by $\L^n(X)$ ($\DC^n(X)$), respectively.
 \end{enumerate} 
\end{definition}

\begin{lemma}\label{zaklvl}
Let $X$ be a Banach space, $F \subset X$  a space of dimension $n$ ($1\leq n < \dim X$), and $A \subset X$. Then the following properties
 are equivalent.
 \begin{enumerate}
 \item[(i)]
 $A$ is an $F$-Lipschitz surface (resp. an $F$-d.c. surface, an $F$-Lipschitz d.c. surface,  an $F$-Lipschitz locally d.c. surface).
 \item[(ii)]
 There exists a topological complement $\wt E$ of $F$ such that 
 $\wt \pi|_A : A \to \wt E$ is a bijection  and  $(\wt \pi|_A)^{-1}$
  is Lipschitz (resp. d.c., etc.), where $\wt \pi :=\pi_{\wt E, F}$.
 \item[(iii)] 
 If $X = F \oplus E$ and  $\pi :=\pi_{E,F}$, then $\pi|_A : A \to E$ is a bijection  and  $(\pi|_A)^{-1}$
  is Lipschitz (resp. d.c., etc.).
 \item[(iv)]
 If $X = F \oplus E$, then there exists
  a Lipschitz mapping (resp. a
 d.c. mapping, etc.) $\vf: E \to F$ such that  $A = \{x+\vf(x): \ x \in E\}$.
 \end{enumerate}
 \end{lemma}
\begin{proof}
 In the proof we use Lemma \ref{fakta}(i)-(iii).
 
 If (i) holds, then there exists a topological complement $\wt E$ of $F$ and a Lipschitz (d.c., etc.) mapping
  $\wt \vf: \wt E \to F$ such that  $A = \{x+ \wt \vf(x): \ x \in \wt E\}$. Set 
$\wt \pi :=\pi_{\wt E,F}$. Then clearly $\wt \pi|_A: A \to \wt E$ is a bijection and $(\wt \pi|_A)^{-1}$ is
 Lipschitz (d.c., etc.), since $(\wt \pi|_A)^{-1}(\wt e) = \wt e + \wt \vf(\wt e)$. 
 
 Now let $\wt E$ be as in (ii), and 
 let $E$ and $\pi$ be
  as in (iii). Since $\pi|_{\wt E}: \wt E \to E$ is clearly a linear isomorphism, $(\pi|_{\wt E})^{-1} =
   \wt \pi|_E$, $\pi|_A = (\pi|_{\wt E}) \circ (\wt \pi|_A)$ and $(\pi|_A)^{-1} = (\wt \pi|_A)^{-1} \circ
    (\wt \pi|_E)$, we easily obtain (iii).
 
    Letting $\vf (x):= (\pi|_A)^{-1}(x) - x$ for $x \in E$, we easily see that (iii) implies (iv). The
     implication (iv) $\Rightarrow$ (i) is trivial.
\end{proof}

\begin{remark}\label{odim}
\begin{enumerate}
\item Every Lipschitz surface of codimension $n$ in $X$ is clearly a closed subset of $X$.
\item If $S\subset X$ is a Lipschitz (resp. d.c., etc.) surface of codimension  $n \geq 2$, then $S$
 is a subset of a Lipschitz (resp. d.c., etc.) surface of codimension $n-1$. Indeed, suppose that
  $S = \{x+ \vf(x):\ x \in E\}$, where $\vf: E \to F$, $X = E \oplus F$, and $F$ is $n$-dimensional.
  Choose $0 \neq v \in F$, and write $F = \span\{v\} \oplus \wt{F}$. Set $\wt E := E + \span\{v\}$ and,
  for $x \in \wt E$, define $\wt \vf (x) := \pi_{\wt F, \wt E}\left(\vf(\pi_{E,F}(x))\right)$.
   Set  $\wt S := \{y + \wt \vf (y):\ y \in \wt E\}$. 
 It is easy to see that $S \subset \wt S$ and $\wt \vf: \wt E \to F$ is Lipschitz (resp. d.c., etc) if
  $\vf$ is Lipschitz (resp. d.c., etc.).  
 
 Consequently, if $\dim X >n\geq 2$, then $\L^n(X) \subset \L^{n-1}(X)$. If $X$ is separable, then this inclusion is proper, see Remark \ref{inkl}, which shows that no Lipschitz surface of codimension $n-1$ belongs to $\L^n(X)$
  (if $\dim X < \infty$, it is sufficient to use by the obvious way basic properties of Hausdorff dimension). 
 \item
 If $X$ is separable, then the $\sigma$-ideal $\DC^n(X)$ coincides with the $\sigma$-ideal generated by Lipschitz d.c. surfaces (or
  Lipschitz locally d.c. surfaces, or locally d.c. surfaces).
   It easily follows from local Lipschitzness of d.c. functions, from the well-known fact that each Lipschitz convex function defined
    an an open ball in a space $E$ can be extended to a Lipschitz convex function on $E$, and from separability of $X$.
 \item 
  It is not difficult to show that  $\DC^n(X)$ is a proper subset of $\L^n(X)$ (if $\dim X > n \geq 1$); see \cite[p. 295]{Za78a}
   for $n=1$.           
\end{enumerate} 

\end{remark}

\begin{remark}\label{closed}
 Suppose that $X = E \oplus F$ and $F$ is finite dimensional.
An easy argument using local compactness of $F$
 shows that $\pi_{E,F}(A)$ is closed in $E$ whenever $A$ 
  is closed and bounded in $X$. Consequently,
  $\pi_{E,F}(A)$ is an $F_{\sigma}$ subset of $E$ whenever $A$ 
  is closed  in $X$.
  \end{remark}

We will need the following well-known easy consequence of the Brouwer's
 Invariance of Domain Theorem. Because of the lack of a suitable reference,
  we present a short proof.
  
\begin{lemma}\label{bro}
Let $C$, $\wt C$ be   Banach spaces with $0 < \dim C = \dim \wt C < \infty$ and let $f: \wt C \to C$
 be an injective continuous mapping such that $f^{-1}: f(\wt C) 
  \to \wt C$ is Lipschitz. Then $f(\wt C) = C$.
\end{lemma}
\begin{proof}
We can clearly suppose that $C= \wt C$ and $X:= C = \wt C$ is an Euclidean space.
The
 Brouwer's
 Invariance of Domain Theorem implies that $f(X)$ is open
  in $X$. Let $y_n \to y$, where $y_n \in f(X)$. Then $(y_n)$
   is bounded and, since $f^{-1}$ is Lipschitz, 
    $(x_n) := (f^{-1}(y_n))$ is bounded as well.
  Choose a subsequence  $x_{n_k} \to x \in X$. Then 
   $f(x_{n_k}) = y_{n_k} \to f(x) = y$. Thus, we have proved that
    $f(X)   $ is closed; the connectivity of $X$ implies
     $f(X) = X$.                                       
\end{proof}

\begin{proposition}\label{brower}
Let $X$ be a Banach space, $S \subset X$  a Lipschitz surface of codimension $n$, and let  $X =  D \oplus F$
 with  $\dim F = n$. Let  $\psi = \pi_{D,F}|_S : S \to D$ be injective and $\psi^{-1}: \psi(S) \to S$
  be Lipschitz. Then $S$ is an $F$-Lipschitz surface.
  Moreover, if $S$ is a  Lipschitz locally d.c. surface of codimension $n$, then 
  $S$ is an $F$-Lipschitz locally d.c. surface.
\end{proposition}
\begin{proof}
Choose an $n$-dimensional space $\wt F$ such that $S$ is an $\wt F$-Lipschitz surface. Since the case $F = \wt F$
 is obvious by Lemma \ref{zaklvl}, we suppose $F \neq \wt F$. Put $K: = F \cap \wt F$ and choose spaces $C$, $\wt C$
  such that $F = K \oplus \wt C$ and $\wt F = K \oplus C$. Then clearly  $1 \leq \dim C = \dim \wt C < \infty$. Choose a topological complement $Z$ of the (finite dimensional) space $F + \wt F = K \oplus C  \oplus \wt C$ and denote
   $E := Z \oplus C$, $\wt E := Z \oplus \wt C$. Clearly $X = F \oplus E = \wt F \oplus \wt E$.

By Lemma \ref{zaklvl}, $\wt \vf := \pi_{\wt E,\wt F}|_S : S \to \wt E$
  is a bilipschitz bijection. It is easy to see (proceeding similarly as in the proof of Lemma \ref{zaklvl})
that  $ \vf := \pi_{ E, F}|_S : S \to  E$ is injective and $\vf^{-1}: \vf(S) \to S$
  is Lipschitz. So Lemma \ref{zaklvl} implies that, to prove the first part of the assertion, it is sufficient
   to verify $\vf (S) =E$. 
   
   To this end 
choose an arbitrary  $e \in E$ and write  $e = z + c$, where $z \in Z$ and $c \in C$.
     For each $x \in \wt C$, put  $f(x) := \vf \circ (\wt \vf)^{-1} (x+z) - z$. Clearly $f(x) \in (F + \wt F) \cap E = C$; so
   $f : \wt C \to C$. It  is easy to see that $f$ is continuous injective  and $f^{-1}(y) = \wt \vf \circ \vf^{-1} (y+z) - z$ for each
      $y \in f(\wt C)$. Consequently, $f^{-1}$ is Lipschitz, and so $f(\wt C)=C$ by Lemma \ref{bro}. For
       $\wt c := f^{-1}(c)$ we have  $\vf ((\wt \vf)^{-1}(\wt c+z)) = c+z = e$; so $\vf(S)=E$. 

To prove the second part of the assertion, we  suppose that $(\wt \vf)^{-1}: \wt E \to X$ is moreover locally d.c.
 Then  $g := \vf \circ (\wt \vf)^{-1} = \pi_{E,F} \circ (\wt \vf)^{-1} $ is clearly Lipschitz and it is locally d.c. by Lemma \ref{fakta} (i). Since  $\vf(S) =E$, we have
  that $g: \wt E \to E$ is a bijection and $g^{-1} = \wt \vf \circ \vf^{-1}$ is Lipschitz. Choose a linear
   bijection $L: \wt C \to C$, and let $G: \wt E \to E$ be the mapping which assigns to a point
    $\wt e = \wt c + z$ ($\wt c \in \wt C$, $z \in Z$) the point $G(\wt e) := L(\wt c) + z$. Then clearly $G$ is a linear
     isomorphism. Since  $G(\wt e)- \wt e \in C+\wt C$ and  $g (\wt e)-\wt e \in F + \wt F$, we obtain that
      $g- G$ has a finite dimensional range. Consequently,  Lemma \ref{fakta}(v) implies that $g^{-1}$
       is locally d.c. Thus, Lemma \ref{fakta}(iv) implies that $\vf^{-1} = (\wt \vf)^{-1} \circ g^{-1}$ is locally
        d.c. So, Lemma \ref{zaklvl} implies that  $S$ is an $F$-Lipschitz locally d.c. surface.

\end{proof}

\begin{proposition}\label{per}
Let $X$ be a  Banach space,
$F\subset X$   an $n$-dimensional space, and $A \subset X$  an
 $F$-Lipschitz (resp. $F$-Lipschitz locally d.c.) surface. Then there exists $\ep >0$
  such that if
 $\wt{F} \subset X$ is an $n$-dimensional space with $\gamma(F,\wt F) < \ep$,
  then   $A$ is an $\wt F$-Lipschitz (resp. $\wt F$-Lipschitz locally d.c.) surface.
\end{proposition}
\begin{proof}
Choose  $E$ such that $X = E \oplus F$, and choose $K \geq 1$ such that the 
canonical mapping $\gamma: E \oplus F \to E \times F$ is 
$K$-bilipschitz.
Choose a corresponding $\omega > 0$ by Lemma \ref{kdk}. 
Denote
$\pi := \pi_{E,F}$, and choose $L \geq 1$ such that $(\pi 
|_A)^{-1}$ is Lipschitz with constant $L$.
Choose
 $\ep  > 0$ such that $\ep < \omega$ and
 \begin{equation}\label{tro}
 2K^2 L \ep < 1/2. 
\end{equation}
Now suppose that  an $n$-dimensional space $\wt F$ with $\gamma(F,\wt F) < \ep$ be given. Since $\ep < \omega$,
 we have that $X = E \oplus \wt F$ 
 and the canonical mapping $\wt \gamma: E \oplus \wt F
  \to E \times \wt F$ is $2K$-bilipschitz.
By Proposition \ref{brower},  it is sufficient to prove that, putting
 $\wt \pi := \pi_{E,\wt F}$, the mapping
  $(\wt \pi|_A)^{-1}$ is Lipschitz with constant $2L$;
 i.e., that
 \begin{equation}\label{zav}
\|x-y\| \leq 2L \|\wt \pi (x) - \wt \pi (y)\| = 2 L \|\wt \pi
 (x-y)\|,\ \ \ x,y \in A. 
\end{equation}  
Thus, consider $x,y \in A,\ x \neq y$, and write $x-y = e_1 + f = e_2 + 
\wt f$, where $e_1= \pi(x-y) \in E$, $e_2 = \wt \pi (x-y) \in E$, $f \in F$, 
and $\wt f \in \wt F$. 
 We know that $\|x-y\| \leq L \|e_1\|$ and so
   $\|\wt f\| \leq 2K \|x-y\| \leq 2KL \|e_1\|$.
  
If $\wt f = 0$, then (\ref{zav}) is obvious. If $\wt f \neq 0$,
put $\wt z := \|\wt f\|^{-1} \wt f$, and find $z \in F$ such
 that $\|\wt z - z\| \leq \ep$. Then $f_2 := \|\wt f\| z$ 
 satisfies  $\|\wt f - f_2\| \leq \ep \|\wt f \|$, and so
 $$ K^{-1} \|e_1-e_2\| \leq \|e_1-e_2+f-f_2\| = \| \wt f - f_2\| \leq
  \ep \|\wt f\| \leq 2KL\ep \|e_1\|.$$
  Thus, by (\ref{tro}), we obtain $\|e_1-e_2\| \leq \|e_1\|/2$, and so
   $\|e_2\|\geq  \|e_1\|/2$. Therefore
   $ \|x-y\|  \leq L \|e_1\| \leq 2L \|e_2\|,$
    which proves (\ref{zav}) and finishes the proof.
\end{proof}

\begin{remark}\label{nevim}
I do not know whether the analogies of Proposition \ref{brower} and Proposition \ref{per} hold for Lipschitz d.c. surfaces.
\end{remark}

\section{Projections of Lipschitz surfaces of finite codimension}

\begin{definition}\label{chv}
Let $X$ be a separable Banach space, and let a finite-dimensional space $V \subset X$ be given. We define the following classes of sets:
\begin{enumerate}
\item[(i)] 
$\cal A(V)$ is the system of all Borel sets $B \subset X$ such that $V \cap (B+a)$ is Lebesgue null (in $V$) for each
 $a \in X$. For $0\neq v \in X$, we put $\cal A(v) := \cal A (\span\{v\})$.
\item[(ii)] 
$\cal A^* (V,\ep)$ (where $0 < \ep < 1$) is the system of all Borel sets $B \subset X$ such that 
$B \in \cal A(W)$ for each space $W$ with $\gamma(V,W) < \ep$, and $\cal A^*(V)$ is the system
 of all sets $B$ such that $B = \bigcup_{k=1}^{\infty} B_k$, where $B_k \in \cal A^*(V, \ep_k)$ for some
  $0 < \ep_k <1$.
\item[(iii)] 
$\cal C^*_d$  (where $d \in \N$) is the system of those $B \subset X$ that can be written as
$B = \bigcup_{k=1}^{\infty} B_k$, where each $B_k$ belongs to $\cal A^* (V_k)$ for some $V_k$ with
 $\dim V_k = d$.
\item[(iv)] 
$\cal A$ is the system of those $B \subset X$ that can be, for every  complete sequence $(v_k)$ in $X$,
 written as  $B = \bigcup_{k=1}^{\infty} B_k$, where each $B_k$ belongs to $\cal A(v_k)$.
\end{enumerate}
\end{definition}
Note that $\cal C^*_1$ coincide with $\cal C^*$ from \cite{PZ} and $\cal A$ is the system of all
 Aronszajn null sets. For basic properties of sets from $\cal A$ see \cite{BL}.
Lemma \ref{gaptop} easily implies the following fact.
\begin{lemma}\label{isom}
Let $X$, $Y$ be Banach spaces and $F: X \to Y$  a linear isomorphism. Let $S\subset X$ belong
 to $\A^*(V,\delta)$. Then there exists $\ep>0$ such that $F(S) \in \A^*(F(V), \ep)$ (in the space $Y$).
\end{lemma}

\begin{proposition}\label{CKA}
Let $X$ be a separable infinite dimensional Banach space and $k \in \N$. 
Then $\cal C_1^* \subset\cal C_2^* \subset\dots \subset \cal A$
 and all inclusions are proper.
\end{proposition}
\begin{proof}
To prove the inclusions $\cal C_d^* \subset \cal A$, it is sufficient to show that 
$\A^*(V,\ep) \subset \A$   whenever $V \subset X$ is a $d$-dimensional space and $\ep>0$.
Let  $V$, $\ep$ and $B \in \A^*(V,\ep)$ be given. Choose a basis $(v_1,\dots,v_d)$ of $V$ and consider an arbitrary
 complete sequence $(u_i)$ in $X$. Choose a $\delta>0$ that corresponds to $(v_1,\dots,v_d)$ and $\ep$
  by Lemma \ref{baze}. We can clearly choose $n \in \N$ and vectors $w_1,\dots w_d$ in $U:= \span \{u_1,\dots,u_n\}$
   such that      $\|w_i - v_i\| < \delta,\ i=1,\dots,d$. Then, denoting  $W := \span \{w_1,\dots,w_d\}$, we
    have  $\gamma(V,W) < \ep$, and so  $B \in \A(W)$. Consequently, by the Fubini theorem, $B \in \A(U)$. Using 
    \cite[Proposition 6.29]{BL}, we easily obtain that $B$ can be decomposed as $B = \bigcup_{i=1}^n B_i$, where
     $B_i \in \A(u_i)$. So, $B \in \A$, and $\cal C_d^* \subset \cal A$ is proved.
     
     To prove $\cal C_d^* \subset \cal C_{d+1}^*$, consider a $B \in \cal A^*(V,\ep)$, where $\dim V = d$ and
      $1>\ep >0$. Choose a basis $v_1,\dots,v_d$ of $V$ with $\|v_i\|=1$ and find a corresponding $\delta >0$ by
       Lemma \ref{baze}. Now choose an arbitrary $Z \supset V$ with $\dim Z = d+1$. To prove $B \in \cal A^*(Z,\delta)$,
        consider an arbitrary $(d+1)$-dimensional  $W$ with $\gamma(W,Z) < \delta$. By the definition of $\gamma$,  find
         $w_1,\dots,w_d \in W$ with $\|w_1-v_1\| < \delta,\dots,\|w_d-v_d\|< \delta$, and set $\wt W := \span\{w_1,\dots,w_d\}$.
         The choice of $\delta$ implies that $\gamma(\wt W,V) < \ep$, and so $\wt W \cap (B+a)$ is Lebesgue null in $\wt W$ for each $a \in X$.
          Consequently, by Fubini theorem, $ W \cap (B+a)$ is Lebesgue null in $ W$ for each $a \in X$. So $B \in \cal A^*(Z,\delta)$, 
           and $\cal C_d^* \subset \cal C_{d+1}^*$ follows.
     
     A construction  
     of a set in $\A \setminus \cal C_1^*$  is presented in the proof of
      \cite[Proposition 13]{PZ}. Moreover, it is shown in \cite{PZ} that this set ($F_2(I)$) meets any $2$-dimensional
       affine space in a $2$-dimensional Lebesgue null set, which shows that even 
        $\cal C_2^* \setminus \cal C_1^* \neq \emptyset$.
        It is not difficult to modify that construction and obtain a set
       in $ \cal C_{d+1}^* \setminus \cal C_d^*$ for each $d$ (see Remark \ref{obec}). However, since the notation is somewhat complicated
        in the general case, we will give a detailed proof for $d=2$ only.
        
       Our construction  starts quite similarly as the construction of a set from
         $\A \setminus \cal C^*$ 
        on p. 20 of \cite{PZ}.
        Namely, by the same procedure as in \cite{PZ} we can define positive numbers
          $c_0, c_1, c_2,\dots$ and nonzero vectors $u_0, u_1,u_2,\dots$ in $X$ such that
           both  $\{u_{6n-3}:\ n \in \N\}$ and  $\{u_{6n}:\ n \in \N\}$ are dense in $X$, and the formula
           $F(x) = \sum_{k=0}^{\infty} c_k x_{k+1} u_k$ (where $x= (x_1, x_2,\dots)$) defines  a linear injective mapping of $\ell_{\infty}$
            to $X$.

        As in \cite{PZ}, we set $I:= \{x \in \ell_{\infty}:\ 1 \leq x_k \leq 2\}$, and equip $I$
             with the topology of pointwise convergence (so it is a compact metrizable space) and with the 
              measure $\mu$ defined as the product of countably many copies of the Lebesgue measure on $[1,2]$.

     Choose two sequences  $\xi^1_1, \xi^1_2,\dots$ and $\xi^2_1, \xi^2_2,\dots$ such that
     $0 < \xi^1_j < 1/(j+1)!$, $0 < \xi^2_j < 1/(j+1)!$, and
     $$ \lim_{k \to \infty} \sum_{j=k}^{\infty}  c_{3j-2} \xi^1_j 2^j \|u_{3j-2}\|/c_{6k-3} =0,\ \
     \lim_{k \to \infty} \sum_{j=k}^{\infty}  c_{3j-1} \xi^2_j 2^j \|u_{3j-1}\|/c_{6k} =0.$$ 
     Now, for $x\in I$, set
     $$ G(x) = \sum_{k=1}^{\infty} c_{3k-2} \xi_k^1 x_1x_3\dots x_{2k-1}u_{3k-2} + 
         \sum_{k=1}^{\infty} c_{3k-1} \xi_k^2 x_2x_4\dots x_{2k}u_{3k-1} + 
     \sum_{k=1}^{\infty} x_k c_{3k} u_{3k}.$$
     Then $G: I \to X$ is a continuous mapping. Indeed, we have $G = F \circ H$, where 
     $$ H(x_1,x_2,\dots) := (0, \xi_1^1 x_1, \xi_1^2 x_2, x_1, \xi_2^1 x_1 x_3, \xi_2^2 x_2 x_4, x_2, 
      \xi_3^1 x_1 x_3 x_5, \xi_3^2 x_2 x_4 x_6, x_3,\dots),$$ 
     and  $H: I \to \ell_{\infty}$ is  clearly continuous. So, $G(I)$ is compact.
     
     Let $e_j$ be the $j$-th member of the canonical basis of $\ell_{\infty}$.
     Observe that if $x\in I$, $k_1,\, k_2 \in \N$, $t,\, \tau      \in \R$ and $x+ t e_{2k_1-1} + \tau e_{2 k_2} \in I$, then
      $G(x+ t e_{2k_1-1} + \tau e_{2 k_2}) = G(x) + t v_{k_1}(x) + \tau w_{k_2} (x)$, where
      $$ v_k(x) := \sum_{j=k}^{\infty} c_{3j-2} \xi_j^1 (x_1 x_3 \dots x_{2j-1}/x_{2k-1}) u_{3j-2} + c_{6k-3} u_{6k-3},$$
       $$  w_k(x) := \sum_{j=k}^{\infty} c_{3j-1} \xi_j^2 (x_2 x_4 \dots x_{2j}/x_{2k}) u_{3j-1} + c_{6k} u_{6k}.$$             
     
     Now consider $x,\, y \in I$ such that $x \neq y$ and $tG(x) + (1-t)G(y) \in G(I)$ for infinitely
      many real $t$. Since $F$ is a linear injection of $\ell_{\infty}$ to $X$, for any such $t$ we
       have $t H(x) + (1-t)H(y) = H(z)$ for some $z \in I$. Considering $(3k+1)$-th coordinates
        of $H(z)$ we obtain $z = tx + (1-t)y$. Consequently, considering $(3k-1)$-th and $3k$-th coordinates of
         $H(z)$, we obtain that, for each $k \in \N$,
         $$ t x_1 x_3\dots x_{2k-1} + (1-t)  y_1 y_3\dots y_{2k-1} = (tx_1 +(1-t)y_1)\cdots (tx_{2k-1} + (1-t) y_{2k-1}), $$           
      $$ t x_2 x_4\dots x_{2k} + (1-t)  y_2 y_4\dots y_{2k} = (tx_2 +(1-t)y_2)\cdots (tx_{2k} + (1-t) y_{2k}). $$        Since the above equalities hold for infinitely many $t$, we infer that $x$ and $y$ differ at most in one
       odd coordinate and at most in one even coordinate (otherwise one of right sides, for sufiicintly large $k$,
        is a polynomial in $t$ of  degree grater than one, which is impossible).
        Consequently, there exist $k_1, k_2 \in \N$ such that  $y \in x + \span \{e_{2k_1-1}, e_{2k_2}\}$;
         so $G(y) \in G(x) + \span \{v_{k_1}(x), w_{k_2}(x)\}$. 
         
         The above analysis shows that the set of lines, which contain any fixed point $G(x),\ x \in I$, and meet the 
          set $G(I)$ in an infinite set, can be covered by countably many planes containing $G(x)$. Therefore
           $G(I)$ meets any $3$-dimensional affine subspace of $X$ in a set of three dimensional Lebesgue measure
            zero. Consequently, $G(I) \in \cal C^*_3$.
            
            Now suppose that  $G(I) \in \cal C^*_2$, hence  $G(I) = \bigcup_{n=1}^{\infty} B_n$, where
             $B_n \in \cal A^*(V_n, \ep_n)$ and $V_n$ are two-dimensional subspaces of $X$. 
              Write $V_n = \span \{p_n, q_n\}$ and choose $\delta_n > 0$ (by Lemma \ref{baze}) such that 
               $\gamma(V_n, \span\{v, w\}) < \ep_n$ whenever $\|v-p_n\| < \delta_n$ and $\|w-q_n\| < \delta_n$.
                For any given $n$ find $k_1, k_2$ such that
               $$  \sum_{j=k_1}^{\infty} 2^j c_{3j-2} \xi^1_j \|u_{3j-2}\| < c_{6k_1-3} \delta_n/2,\ 
   \  \|u_{6k_1-3} - p_n\| < \delta_n/2, $$
    $$  \sum_{j=k_2}^{\infty} 2^j c_{3j-1} \xi^2_j \|u_{3j-1}\| < c_{6k_2} \delta_n/2,\ \ 
    \|u_{6k_2} - q_n\| < \delta_n/2. $$        
         For any $x \in I$ we have
     $$ \|v_{k_1}(x) - c_{6k_1-3} p_n\| \leq  \sum_{j=k_1}^{\infty} 2^j c_{3j-2} \xi^1_j \|u_{3j-2}\|+   c_{6k_1-3}\|u_{6k_1-3}-p_n\| < c_{6k_1-3} \delta_n,$$
     $$ \|w_{k_2}(x) - c_{6k_2} q_n\| \leq  \sum_{j=k_2}^{\infty} 2^j c_{3j-1} \xi^1_j\|u_{3j-1}\| +   c_{6k_2}\|u_{6k_2}-q_n\| < c_{6k_2} \delta_n.$$     
    So  $ \|v_{k_1}(x)/c_{6k_1-3}- p_n\| < \delta_n$ and  $ \|w_{k_2}(x)/c_{6k_2} - q_n\| < \delta_n$, which
     shows that the plane  $G(x) + \span \{v_{k_1}(x), w_{k_2}(x)\}$ meets $B_n$ in a  two-dimensional Lebesgue
      null set. Hence the set
      $$ \{(t,\tau):\ x + t e_{2k_1-1} + \tau e_{2k_2} \in G^{-1} (B_n) \} = \{(t,\tau):\  G(x) + t v_{k_1}(x) + \tau w_{k_2}(x) \in B_n\}$$
      is Lebesgue null, and Fubini theorem gives $\mu(G^{-1}(B_n))=0$. (Note that $G^{-1}(B_n)$ is Borel, since $G$
       is continuous.) But this contradicts  $I = \bigcup_{n=1}^{\infty} G^{-1}(B_n)$, and we infer that
        $G(I) \notin \cal C^*_2$.
         \end{proof}
        \begin{remark}\label{obec}
        For an arbitrary $d \in \N$, we obtain as above that $G_d(I) \in  \cal C^*_{d+1}\setminus  \cal C^*_d$, where $G_d = F \circ H_d$,
        $$ H_d(x) := (0,\xi^1_1 x_1,\dots,\xi^d_1 x_d, x_1, \xi^1_2 x_1 x_{d+1},\dots,\xi^d_2 x_d x_{2d},x_2, \xi^1_3 x_1 x_{d+1} x_{2d+1},\dots),$$
         and  $(\xi^1_i), \dots, (\xi^d_i)$ are suitably chosen sequences.
     \end{remark}

\begin{proposition}\label{proj}
Let $X$ be a separable infinite dimensional space, $S$  a Lipschitz surface of codimension $n \geq 2$, and
$P:X \to Y$  a continuous linear mapping onto a Banach space $Y$ such that $\dim(\ker(P)) < n$.
 Then there exists an $n$-dimensional  space $D \subset Y$ and $0 < \ep <1$ such
  that $P(S) \in \cal C^*(D, \ep)$ in $Y$. Consequently,  $P(S)$ is a first category subset of $Y$ which is
    Aronszajn null  in $Y$.
\end{proposition}
\begin{proof}
 Denote  $K:= \ker P$. Choose a space $F\subset X$ such that $\dim F = n$ and $S$ is an $F$-Lipschitz
  surface.  Using Lemma \ref{per},  Lemma \ref{zname}
  and Lemma \ref{prpr}, we can choose  a space $V$ with $\dim V= n$ such that $S$ is an $V$-Lipschitz surface
   and $V \cap K = \{0\}$. Choose a closed space $H \subset X$ such that $X = H \oplus (K \oplus V)$. Denoting
    $Z :  = H \oplus V$, we have $X = Z \oplus K$. Set $\pi:= \pi_{Z,K}$. Using Lemma \ref{per} and  Lemma \ref{zname}, we find 
     $0<\delta <1$ such that $\gamma(V,W) < \delta$ implies that $S$ is a $W$-Lipschitz surface and
      $W \oplus (H \oplus K) =  X$. Now consider on arbitrary $W \subset Z$ such that $\gamma(V,W) < \delta$. 
 We can choose a Lipschitz mapping $\vf: H \oplus K \to W $ such that $S= \{ h+k + \vf(h+k):\ h \in H,\,  k \in K\}$.
  Consequently, $\pi(S)= \{ h + \vf(h+k):\ h \in H,\,  k \in K\}$. Now consider an arbitrary
   $a= h_0 + w_0 \in Z$. Then 
   $(\pi(S)+ a) \cap W= \{ w_0 + \vf(-h_0+k):\  k \in K\}$. Since the mapping $\psi: K \to W$ defined by
    $\psi(k):=  w_0 + \vf(-h_0+k)$ is Lipschitz and $\dim K < \dim W$, we obtain that                  $(\pi(S)+ a) \cap W$ is Lebesgue null in $W$. Since $\pi(S)$ is an $F_{\sigma}$ set by Remark \ref{odim}(i) and Remark \ref{closed}, we obtain   that $\pi(S) \in \cal C^*(V, \delta)$ in $Z$.
    Since $F:= P|_Z$ is a linear isomorphism with $F(\pi(S)) = P(S)$, Lemma \ref{isom} implies that
     $P(S) \in \cal C^*(D, \ep)$ for $D:=F(V)$ and some $\ep>0$. Consequently, $P(S)$ is  Aronszajn null in $Y$ by Lemma \ref{CKA}. Thus, $\vn P(S) = \emptyset$. Since $P(S)$ is an $F_{\sigma}$ set, we obtain that 
      $P(S)$ is a first category set.
     \end{proof}

As an immediate consequence, we obtain the following result.
\begin{proposition}\label{proj2}
Let $X$ be a separable infinite dimensional Banach space, $n \geq 2$, $A \in \L^n(X)$, and let
$P:X \to Y$ be a continuous linear mapping onto a Banach space $Y$ such that $\dim(\ker(P)) < n$.
 Then $P(S)$ is a subset of a set from $\cal C^*_n$ in $Y$. Consequently,  $P(S)$ is a first category subset of $Y$ which is
   a subset of an Aronszajn null set in $Y$.
\end{proposition}

\begin{remark}\label{DH}
Let $X$, $Y$, $P$ and $n$ be as in Proposition \ref{proj2}.
\begin{enumerate}
\item 
 Let $f$ be a continuous convex function on $X$ and $B_n := \{x\in X:\ \dim (\partial f(x)) \geq n\}$.
 Then \cite[Theorem 1.3.]{PreH} states that $P(A)$ is a first category set. 
Using results of \cite{Za79}, it is easy to see that \cite[Theorem 1.3.]{PreH} is  equivalent
 to the statement that $P(A)$ is a first category for each $A \in \DC^n(X)$, but the proof of \cite{PreH} is direct, it does not use \cite{Za79}.
 \item  The result that $P(A)$ is a first category for each $A \in \L^n(X)$ is due to Heisler \cite{HeD}.
 \item  An example from \cite{HeD} shows that there exists $A \in \DC^n(X)$ such that $P(A) \notin \L^1(Y)$.
 \item It is not known whether $P(A)$ is $\sigma$-porous or $\Gamma$-null for
  each  $A \in \L^n(X)$ (resp. $A \in \DC^n(X)$). The negative answer seems to be probable.
\end{enumerate}
\end{remark}

\begin{remark}\label{inkl}
Let $X$ be a separable infinite dimensional space. Proposition \ref{proj2} easily implies that the inclusions
  $\L^n(X) \subset \L^{n-1}(X)$ ($n>1$) are
proper. Indeed, no  Lipschitz surface $S$ of codimension $n-1$ can belong to $\L^n(X)$, since there is a surjective continuous linear 
 projection of $S$ on a space $E$ of codimension $n-1$.
 \end{remark}

 Proposition \ref{proj2}  implies the following result
  which improves   both \cite[Theorem 1.3.]{PreH} and \cite[Theorem 5.6.]{HeD}.

  \begin{theorem}\label{mon}
Let $X$ be a separable infinite dimensional space, $n \geq 2$, and let $T: X \to X^*$ be a monotone
 (mutivalued) operator. Denote by $B_n$ the set of all $x\in X$ for which the convex cover of $T(x)$
  is at least $n$-dimensional. Let
$P:X \to Y$ be a continuous linear mapping onto a Banach space $Y$ such that $\dim(\ker(P)) < n$.
 Then $P(B_n)$ is a subset of a set from $\cal C^*_n$ in $Y$. Consequently,  $P(B_n)$ is a first category subset of $Y$ which is
   a subset of an Aronszajn null set in $Y$. 
  \end{theorem}
  \begin{proof}
  Since $B_n \in \L^n(X)$ by \cite{Za78b}, the assertion follows from Proposition \ref{proj2}.
  
  \end{proof}

        \bigskip
        
        {\bf Acknowledgements.}\ \ The research was supported by the institutional
         grant MSM 0021620839 and by the grant GA\v CR 201/06/0198.

\bigskip
\bigskip  
{\it Authors adress:\ \ Lud\v ek Zaj\'\i\v cek}, Charles University,
Faculty of Mathematics and Physics, 
Sokolovsk\'a 83,
186 75 Praha 8, 
Czech Republic,
e-mail:\ zajicek@karlin.mff.cuni.cz.

\end{document}